\documentclass[11pt]{article}

\usepackage{amsmath,amssymb,amsthm,amsfonts}

\numberwithin{equation}{section}

\usepackage{mathtools}
\usepackage{mathrsfs}
\usepackage{bm}
\usepackage{graphicx}
\usepackage{float}
\usepackage{multirow}
\usepackage{booktabs}
\usepackage{longtable}
\usepackage{arydshln}
\usepackage{enumerate}
\usepackage{cite}
\usepackage{microtype}
\usepackage{enumitem}
\usepackage{hyperref}
\usepackage[a4paper,margin=2.5cm]{geometry}

\hypersetup{
 colorlinks=true,
 linkcolor=blue,
 citecolor=blue,
 urlcolor=blue,
pdftitle={Finite Sums of Weighted
Composition Operators on the Fock Space},
pdfauthor={Zhaopeng Lin, Rong Yang, Chao Zu, Yufeng Lu},
 pdfsubject={Finite weighted composition sums on the classical Hilbert Fock space}
}
\allowdisplaybreaks[2]
\newcommand{\doi}[1]{\href{https://doi.org/#1}{\nolinkurl{doi:#1}}}

\newtheorem{theorem}{Theorem}[section]

\newtheorem{lemma}[theorem]{Lemma}
\newtheorem{corollary}[theorem]{Corollary}

\newcommand{\C}{\mathbb C}
\newcommand{\F}{\mathcal F}
\newcommand{\Sp}{\mathcal S}

\newcommand{\ip}[2]{\langle #1,#2\rangle}
\newcommand{\norm}[1]{\left\lVert #1\right\rVert}

\title{\Large\bfseries Finite Sums  of Weighted Composition Operators on the Fock Space}
\author{
Zhaopeng Lin,
Yufeng Lu,
Rong Yang\thanks{Corresponding author.},
Chao Zu
}
\date{}

\begin{document}
% Article-wide rule: no tables or decorative frames around formulas/results.
% Existing label identifiers and the established numbering mechanism are preserved.

\maketitle

\begin{center}
\begin{minipage}{14cm}\small
\noindent\textbf{Abstract}\quad
We study finite sums of weighted
composition operators on the classical Fock space.
We show that no cancellation phenomenon occurs
when the composition symbols are pairwise distinct:
such a sum is bounded, compact, or belongs to the
Schatten class $\mathcal S_p$, $0<p<\infty$,
if and only if each summand has the corresponding
property. We also obtain two-sided essential norm
estimates and quantitative Schatten
norm comparisons for fixed distinct
nonconstant symbols. As an application, we obtain corresponding
boundedness, compactness, and Schatten-class
results for finite sums of generalized
Volterra-type integral operators.
\medskip


\end{minipage}
\end{center}

\section{Introduction}\label{sec:introduction}

Composition operators form a fundamental class of
linear operators on spaces of analytic functions. Given an analytic self-map $\varphi$ of a domain, the associated composition operator
is defined by $C_\varphi f=f\circ\varphi$.  
A central problem is to understand how
operator theoretic properties, such as boundedness,
compactness, and essential norms, are determined
by the analytic and geometric properties of the
inducing symbol. 
This problem has been extensively investigated
on Hardy, Bergman, and other classical spaces
of holomorphic functions; see 
\cite{CowenMacCluer1995,Shapiro1993,MacCluerShapiro1986,Shapiro1987}. 

Beyond individual operators, the study of
differences and finite sums of 
composition operators is closely connected
with the topology of the space of composition
operators. Berkson~\cite{Berkson1981}
showed that certain composition operators
on the Hardy space are isolated in the
operator-norm topology.
Shapiro and Sundberg~\cite{ShapiroSundberg1990}
further investigated isolation phenomena
and raised questions concerning compact
differences and connected components. 
These investigations led to the study
of compact differences. In both Hardy and
  Bergman spaces, cancellation may occur:
 there exist distinct symbols $\varphi$ and $\psi$
such that neither $C_\varphi$ nor $C_\psi$
is compact, whereas $C_\varphi-C_\psi$
is compact; see
\cite{ShapiroSundberg1990,Moorhouse2005,
KrieteMoorhouse2007,ChoeChoiKooYang2020}.
Thus, compactness of a difference need not
imply compactness of its individual summands.

The study of compact differences naturally
leads to the more general problem of linear
relations among composition operators
modulo compact operators. Although the
collection of composition operators is not
closed under addition, nontrivial linear
relations may arise in the Calkin algebra.
Moorhouse~\cite{Moorhouse2005}
characterized compact differences on weighted
Bergman spaces and investigated additive
relations modulo compact operators.
Kriete and Moorhouse~\cite{KrieteMoorhouse2007}
further studied finite linear relations
in the Calkin algebra on Hardy and weighted
Bergman spaces. Related results on finite
linear combinations were obtained on Bloch
spaces~\cite{HosokawaNieminenOhno2011}.
These developments raise the question of
how cancellation among several composition
operators affects the compactness of their 
finite sums.

The situation on Fock spaces is substantially
different. Let $\mathcal H(\C)$ denote the
space of entire functions on $\C$.
For $\alpha>0$ and $0<p<\infty$, the Fock
space $\F_\alpha^p$ consists of all
$f\in\mathcal H(\C)$ such that
\begin{equation} 
\|f\|_{p,\alpha}^p
=\frac{\alpha p}{2\pi}
\int_\C |f(z)|^p
e^{-\frac{\alpha p}{2}|z|^2}\,dA(z)
<\infty,\nonumber
\end{equation}
where $dA$ denotes the planar Lebesgue measure.
In particular, $\F_\alpha^2$ is a Hilbert space,
and we write $\|f\|_\alpha=\|f\|_{2,\alpha}$.

In contrast to the Hardy and Bergman spaces
on the unit disk, where every analytic
self-map induces a bounded composition operator,
boundedness on Fock spaces imposes strong
restrictions on the inducing symbol. 
Indeed, Carswell, MacCluer and Schuster
\cite{CarswellMacCluerSchuster2003}
showed that boundedness forces the composition
symbol to be affine.
 More precisely, $C_\varphi$ is
bounded on $\F_\alpha^2$ if and only if
\[
\varphi(z)=az+b,\qquad |a|\le1,
\]
with $b=0$ whenever $|a|=1$; moreover,
$C_\varphi$ is compact if and only if $|a|<1$.

A natural extension of composition operators
is provided by weighted composition operators.
Given entire functions $u$ and $\varphi$, the
weighted composition operator is defined by
\[
W_{u,\varphi}f=u\cdot(f\circ\varphi).
\] 
This class includes composition and
multiplication operators as special cases. 
Weighted composition operators arose naturally
in the study of isometries on Hardy spaces
and commutants of analytic Toeplitz operators;
see \cite{Forelli1964,Cowen1978,Le2014}.
On Fock spaces, the Weyl unitary operators
provide important examples of weighted
composition operators \cite{Le2014}.

For individual weighted composition operators
on the classical Fock space, explicit criteria
for boundedness and compactness are well known. 
Ueki \cite{Ueki2007} obtained characterizations
in terms of Berezin-type integral transforms.
Le \cite{Le2014} subsequently established
more explicit criteria involving the growth
of the weight and the composition symbol,
and investigated normal and isometric weighted
composition operators.
Further developments include operators
between different Fock spaces, essential norm
estimates, Schatten-class membership, and
topological properties; see, for example,
\cite{TienKhoiWeighted2019,Mengestie2013}.
These results provide a detailed understanding
of individual weighted composition operators.

Choe, Izuchi and Koo~\cite{ChoeIzuchiKoo2010}
proved that, for distinct composition symbols
$\varphi$ and $\psi$ and nonzero constants
$a,b$, the sum $aC_\varphi+bC_\psi$ is
bounded (resp.\ compact) if and only if
both summands are bounded (resp.\ compact).
They also obtained the corresponding
Schatten-class results. 
Cho, Choe and Koo~\cite{ChoChoeKoo2014}
extended these results to arbitrary finite
linear combinations
$\sum_{j=1}^N a_jC_{\varphi_j}$
on Fock--Sobolev spaces.
Their proof uses a lower estimate for
finite sums of reproducing kernels with
fixed nonzero coefficients, obtained through
induction and geometric separation arguments.
These results, however, concern constant
coefficients and do not directly apply to
finite sums with arbitrary entire weights.

For weighted composition operators,
Tien and Khoi~\cite{TienKhoiDifference2019}
proved that, for $0<p\le q<\infty$ and
distinct composition symbols, the difference
$W_{u_1,\varphi_1}-W_{u_2,\varphi_2}$
from $\F_1^p$ to $\F_1^q$ is bounded
(resp.\ compact) if and only if both
summands are bounded (resp.\ compact).
They also obtained two-sided essential norm
estimates for bounded differences with
nonconstant composition symbols when
$1<p\le q<\infty$.
Yang~\cite{Yang2024} subsequently established
analogous boundedness and compactness
results for $0<q<p<\infty$. 
However, these results concern only differences
of two weighted composition operators.
The approach of Tien and Khoi relies on
a lower estimate for a linear combination
of two reproducing kernels, obtained by
evaluating at their centers.
As noted in \cite[Remark~2.3]{TienKhoiDifference2019},
this argument does not directly extend
to sums of three or more weighted
composition operators.

To address this gap in the existing literature, we investigate
finite sums of weighted composition operators
on the classical Hilbert Fock space
$\F_\alpha^2$.
 Let $N\ge1$, and let
$u_j,\varphi_j$ be entire functions for
$1\le j\le N$. We consider the formal sum
\begin{equation}\label{eq:formal-sum}
Tf=\sum_{j=1}^N W_{u_j,\varphi_j}f
=\sum_{j=1}^N u_j(f\circ\varphi_j).
\end{equation} 
Since terms with identical composition
symbols can be combined, we may assume
that the composition symbols are pairwise
distinct.

The following are our main results. 
\begin{theorem}\label{T1.1}
Let $T$ be given by \eqref{eq:formal-sum}, and suppose that the composition
symbols $\varphi_1,\ldots,\varphi_N$ are pairwise distinct entire
functions. Then $T$ is bounded on $\F^2_\alpha$ if and only if each
$W_{u_j,\varphi_j}$ is bounded on $\F_\alpha^2$.
\end{theorem}

We next establish an essential norm estimate
and characterize compactness of finite sums. 
For entire functions $u$ and $\varphi$, set
\begin{equation}\label{eq:m-def}
m_{u,\varphi}(z)
=|u(z)|\exp\left(
\frac{\alpha}{2}
\bigl(|\varphi(z)|^2-|z|^2\bigr)
\right).
\end{equation}

\begin{theorem}\label{T1.2}
Let $T$ be given by \eqref{eq:formal-sum},
where the composition symbols
$\varphi_1,\ldots,\varphi_N$ are pairwise
distinct entire functions. 
If $T$ is bounded on $\F_\alpha^2$,
then there exist constants $c,C>0$, 
depending only on $\alpha$ and
$\varphi_1,\ldots,\varphi_N$, 
 such that
\begin{equation}\label{eq:finite-essential} 
c\limsup_{|z|\to\infty}
\left(
\sum_{j=1}^N m_{u_j,\varphi_j}(z)^2
\right)^{1/2}
 \le \|T\|_e 
 \le C\limsup_{|z|\to\infty}
\left(
\sum_{j=1}^N m_{u_j,\varphi_j}(z)^2
\right)^{1/2}. 
\end{equation}
In particular, $T$ is compact on
$\F_\alpha^2$ if and only if
$W_{u_j,\varphi_j}$ is compact on
$\F_\alpha^2$, for $1\le j\le N$.
\end{theorem}

We also characterize Schatten-class membership
and obtain estimates for the Schatten norms.

\begin{theorem}\label{T1.3}
Let $0<p<\infty$ and $T$ be given by \eqref{eq:formal-sum}, and suppose that the composition
symbols $\varphi_1,\ldots,\varphi_N$ are pairwise distinct entire
functions. Then $T\in\Sp_p$ if and only if each
$W_{u_j,\varphi_j}$ belongs to $\Sp_p$.   Moreover, if the symbols are also nonconstant, then there are $c,C>0$, depending
only on $\alpha,p,N$ and the fixed symbols, such that
\begin{equation}\label{eq:finite-schatten-quantitative}
 c\sum_{j=1}^N\norm{W_{u_j,\varphi_j}}_{\Sp_p}^p
 \le\norm T_{\Sp_p}^p
 \le C\sum_{j=1}^N\norm{W_{u_j,\varphi_j}}_{\Sp_p}^p.
\end{equation} 
\end{theorem}

Our approach is based on polynomial
interpolation, which provides a uniform
coefficient estimate for arbitrary finite
sums of reproducing kernels.
Combined with a variance identity and
the growth rigidity of entire functions,
this estimate shows that the pairwise
differences of distinct composition
symbols are affine and uniformly separated
outside a sufficiently large disk. 
We use this separation to show that
boundedness of the sum implies boundedness
of each summand, and the same holds
for compactness. 
For Schatten-class membership, we further
combine uniform estimates for families
of reproducing kernels with a discrete
operator factorization to obtain
quantitative Schatten norm comparisons
for all $0<p<\infty$.
The same approach also applies to finite
sums of generalized Volterra-type
integral operators.

For entire functions $g$ and $\phi$, consider
the generalized Volterra-type integral operator
\begin{equation} 
V_{(g,\phi)}f(z)
=\int_0^z f(\phi(\zeta))g'(\zeta)\,d\zeta. \nonumber
\end{equation}
Mengestie~\cite{Mengestie2014} studied
the boundedness, compactness, and
Schatten-class membership of individual
operators. Mengestie and
Worku~\cite{MengestieWorku2018}
further investigated Schatten-class properties
and differences of two such operators.
More recently, Hao and Tien~\cite{HaoTien2025}
characterized the boundedness and compactness
of differences between Fock spaces. 
Here we obtain corresponding results for
finite sums on $\F_\alpha^2$.

\begin{theorem}\label{T1.4}
Let $N\ge1$ and $g_j,\phi_j\in\mathcal H(\C)$,
$1\le j\le N$, with
$\phi_1,\ldots,\phi_N$ pairwise distinct.
Then, for each of the properties of boundedness,
compactness, and membership in $\Sp_p$
($0<p<\infty$), the sum
\[
\mathcal V=\sum_{j=1}^N V_{(g_j,\phi_j)}
\]
has that property on $\F_\alpha^2$ if and only if
each $V_{(g_j,\phi_j)}$ does.
\end{theorem}

The present paper and the recent work of Tong and Yang
\cite{TongYang2026FiniteSums} were developed independently.
Tong and Yang establish boundedness and compactness
criteria for finite sums of weighted composition
operators between Fock spaces for the full range
$0<p,q<\infty$. The overlap between the two papers
concerns the boundedness and compactness rigidity
results on the Fock space $\F_\alpha^2$. 
The directions of the two works are otherwise complementary: Tong and Yang focus on boundedness and compactness of operators between Fock spaces, whereas the present
paper establishes two-sided essential norm estimates,
Schatten-class characterizations  for all $0<p<\infty$,
as well as corresponding results for finite sums
of generalized Volterra-type integral operators.

The paper is organized as follows.
Section~\ref{sec:preliminaries} recalls
the necessary preliminaries.
Section~\ref{sec:boundedness} establishes
the boundedness result for finite sums.
Section~\ref{sec:compactness} studies
compactness and essential norm estimates.
Section~\ref{sec:finite-schatten}
investigates Schatten-class membership
and quantitative norm comparisons.
Finally, Section~\ref{sec:volterra}
extends corresponding results for finite sums to generalized Volterra-type
operators.

\section{Preliminaries}\label{sec:preliminaries}

\subsection{Fock spaces and weighted composition operators}

For $\alpha>0$, the inner product on
$\F_\alpha^2$ is given by 
\[
 \ip fg=\frac{\alpha}{\pi}\int_\C
 f(z)\overline{g(z)}e^{-\alpha|z|^2}\,dA(z).
\]
The monomials $e_n(z)=\alpha^{n/2}z^n/\sqrt{n!}$,
$n\geq 0 $, form an orthonormal basis. The reproducing kernel
and its normalization are
\begin{equation} 
 K_w(z)=e^{\alpha z\overline w},\qquad
 k_w(z)=e^{\alpha z\overline w-\alpha|w|^2/2},
 \qquad z,w\in\C.  \nonumber 
\end{equation} 
Thus $\|k_w\|_\alpha=1$, and $k_w$
converges weakly to $0$ in $\F_\alpha^2$
as $|w|\to\infty$. 
Moreover,
\begin{equation}\label{eq:kernel-overlap}
 |\ip{k_z}{k_w}|=e^{-\alpha|z-w|^2/2}.
\end{equation}

Whenever $W_{u,\varphi}$ is bounded, the reproducing identity gives
\begin{equation}\label{eq:adjoint-kernel}
 W_{u,\varphi}^*k_z
 =\overline{u(z)}e^{\frac{\alpha}{2}
 (|\varphi(z)|^2-|z|^2)}k_{\varphi(z)},
\end{equation}
so $m_{u,\varphi}(z)=\norm{W_{u,\varphi}^*k_z}_\alpha$.

Let $u,\varphi\in\mathcal H(\C)$, with
$\varphi$ nonconstant. 
By the same argument as in
\cite[Proposition~2.1 and Theorem~2.2]{Le2014},
$W_{u,\varphi}$ is bounded on $\F_\alpha^2$
if and only if
$m_{u,\varphi}\in L^\infty(\C)$. 
Moreover, if $u\not\equiv0$ and
$W_{u,\varphi}$ is bounded, then 
\[
\varphi(z)=az+b,\qquad 0<|a|\le1,
\] for some $a,b\in\C$. 
For a constant symbol $\varphi\equiv b$, we have
\begin{equation} 
W_{u,b}f=uf(b)=\ip{f}{K_b}\,u.
\nonumber 
\end{equation}
Thus, $W_{u,b}$ is bounded on $\F_\alpha^2$
if and only if $u\in\F_\alpha^2$,
in which case it has rank at most one.

For $r>0$, write $
D(w,r)=\{z\in\C:|z-w|<r\}$. 
A set $Z\subset\C$ is $\delta$-separated if
$|z-w|\ge\delta$ for all distinct $z,w\in Z$.
We call $Z$ relatively separated if
\[
\sup_{w\in\C}\#(Z\cap D(w,1))<\infty.
\]

For a countable set $J$ and $0<p<\infty$,
let $\ell^p(J)$ denote the space of all
$p$-summable complex sequences indexed by $J$.

By the $p=2$ case of
\cite[Lemma~4.9]{Zhu2012}
and Hilbert space duality, for every
$\delta>0$ there exists
$C_{\alpha,\delta}>0$ such that 
\[
\left\|\sum_{j\in J}c_jk_{w_j}\right\|_\alpha^2
\le C_{\alpha,\delta}
\sum_{j\in J}|c_j|^2
\]
whenever $(w_j)_{j\in J}$ is
$\delta$-separated and $(c_j)$ is
finitely supported.

\subsection{Schatten ideals}

Let $H_1$ and $H_2$ be complex Hilbert spaces.
For $0<p<\infty$, the Schatten class
$\Sp_p(H_1,H_2)$ consists of all compact
operators $S:H_1\to H_2$ such that
\[
\|S\|_{\Sp_p}^p
=\sum_{n=1}^\infty s_n(S)^p<\infty,
\]
where $s_n(S)$ denotes the $n$-th singular
value of $S$, counted with multiplicity.
When the underlying spaces are clear,
we simply write $\Sp_p$.

For $1\le p<\infty$, $\|\cdot\|_{\Sp_p}$
is a norm, whereas for $0<p<1$ it is
a quasi-norm satisfying
\[
\|S+T\|_{\Sp_p}^p
\le \|S\|_{\Sp_p}^p+\|T\|_{\Sp_p}^p.
\]

The Schatten classes satisfy the
two-sided ideal property. More precisely,
if $A$ and $B$ are bounded operators
such that $ASB$ is well defined,
and $S\in\Sp_p$, then $ASB\in\Sp_p$ and 
\begin{equation}\label{eq:ideal-property}
\|ASB\|_{\Sp_p}
\le \|A\|\,\|S\|_{\Sp_p}\,\|B\|.
\end{equation}
Here $\|\cdot\|$ denotes the operator norm.
In particular, $\|S\|\le\|S\|_{\Sp_p}$.

Moreover, for any finite family
$S_1,\ldots,S_M\in\Sp_p$, we have
\[
\left\|\sum_{j=1}^M S_j\right\|_{\Sp_p}^p
\le M^{\max\{p-1,0\}}
\sum_{j=1}^M\|S_j\|_{\Sp_p}^p.
\]
These properties hold for all $0<p<\infty$;
see \cite[Section~6.3.1]{Pietsch2007}
and \cite[Chapter~2]{Simon2005}.

\section{Boundedness}\label{sec:boundedness}

\begin{lemma}\label{lem:kernel-extraction}
Let $M\ge1$ be an integer, and let
$w_1,\ldots,w_M,c_1,\ldots,c_M\in\C$.
There exists a constant $C_{M,\alpha}\ge1$
such that, for every $1\le j\le M$,
\begin{equation}\label{eq:kernel-extraction}
\left\|\sum_{\ell=1}^M c_\ell K_{w_\ell}\right\|_\alpha
\ge C_{M,\alpha}^{-1}|c_j|e^{\alpha|w_j|^2/2}
\prod_{\ell\ne j}
\frac{|w_j-w_\ell|}{1+|w_j-w_\ell|}.
\end{equation}
\end{lemma}

\begin{proof}
Fix $j$ and assume first that $w_j\ne w_\ell$ for all $\ell\ne j$.  Define
\[
 F_j(\zeta)=k_{w_j}(\zeta)
 \prod_{\ell\ne j}\frac{\zeta-w_\ell}{w_j-w_\ell}.
\]
Then $F_j(w_\ell)=0$ for $\ell\ne j$, and $F_j(w_j)=e^{\alpha|w_j|^2/2}$.  Since the inner product is linear in the first variable, we have
\begin{equation}\label{e2.1}
\left|
\left\langle
F_j,\sum_{\ell=1}^M c_\ell K_{w_\ell}
\right\rangle
\right|
=|c_j|e^{\alpha|w_j|^2/2}.
\end{equation}
Moreover, after the change of variables $\xi=\zeta-w_j$,
\begin{align}
 \norm{F_j}_\alpha^2
 &=\frac{\alpha}{\pi}\int_\C e^{-\alpha|\xi|^2}
 \prod_{\ell\ne j}
 \left|1+\frac{\xi}{w_j-w_\ell}\right|^2\,dA(\xi)\nonumber\\
 &\le C_{M,\alpha}^2
 \prod_{\ell\ne j}
 \left(1+\frac1{|w_j-w_\ell|}\right)^2. \label{e2.2}
\end{align} 
The last inequality follows from
\[
\left|1+\frac{\xi}{d}\right|
\le (1+|\xi|)\left(1+\frac1{|d|}\right),
\qquad d\ne0,
\]
and the integrability of
$(1+|\xi|)^{2(M-1)}e^{-\alpha|\xi|^2}$
over $\C$.
  Combining \eqref{e2.1} and \eqref{e2.2}, Cauchy--Schwarz inequality now gives \eqref{eq:kernel-extraction}.  If one of the points coincides with $w_j$, the asserted lower bound is trivial.
\end{proof}

We also use the following immediate growth consequence of
Le~\cite[Proposition~2.1]{Le2014}. If $F,H\in \mathcal H(\C)$,
$F\not\equiv0$, and
\begin{equation} 
 |F(z)|e^{A|H(z)|^2-B|z|^2}\le C\qquad(z\in\C) \nonumber 
\end{equation}
for some $A>0$, $B\ge0$, and $C>0$, then $H$ is affine.

\begin{proof}[Proof of Theorem \ref{T1.1}]
The sufficiency is immediate. Now, we assume that $T$ is bounded.  If all weights vanish there is nothing to prove. Otherwise, delete zero weights and relabel the remaining terms, so that $u_j\not\equiv0$ for all $j$.  If $N=1$, it is obvious.  We  assume $N\ge2$.

For each $z\in\C$, we have
\begin{equation}\label{eq:Tstar-unnormalized}
 T^*K_z=\sum_{j=1}^N\overline{u_j(z)}K_{\varphi_j(z)},
 \qquad
 \norm{T^*K_z}_\alpha\le\norm{T}e^{\alpha|z|^2/2}.
\end{equation}

Set
\[
 h_{j\ell}=\varphi_j-\varphi_\ell,
 \qquad
 \rho(t)=\frac{t}{1+t}\quad(t\ge0).
\]
Applying Lemma~\ref{lem:kernel-extraction} to \eqref{eq:Tstar-unnormalized} and dividing by $e^{\alpha|z|^2/2}$ yields, for every $j$,
\begin{equation}\label{eq:amplitude-with-separation}
 m_{u_j,\varphi_j}(z)
 \prod_{\ell\ne j}\rho(|h_{j\ell}(z)|)
 \le C_{N,\alpha}\norm{T}.
\end{equation}
Multiplying \eqref{eq:amplitude-with-separation} over $j=1,\ldots,N$ gives
\begin{equation}\label{eq:vandermonde-raw}
 \left(\prod_{j=1}^N m_{u_j,\varphi_j}(z)\right)
 \prod_{j<\ell}\rho(|h_{j\ell}(z)|)^2
 \le C.
\end{equation}
Because the symbols are pairwise distinct, every $h_{j\ell}$ is a nonzero entire function.  Define the nonzero entire function
\[
P(z)
 =\left(\prod_{j=1}^N u_j(z)\right)
 \left(\prod_{j<\ell}h_{j\ell}(z)^2\right).
\]
Expanding the definition of $m_{u_j,\varphi_j}$ in \eqref{eq:vandermonde-raw}, we obtain 
\begin{equation}\label{eq:vandermonde-growth}
 |P(z)|
 \exp\!\left\{\frac{\alpha}{2}
 \left(\sum_{j=1}^N|\varphi_j(z)|^2-N|z|^2\right)\right\}
 \le C\prod_{j<\ell}(1+|h_{j\ell}(z)|)^2.
\end{equation}
For every $\varepsilon>0$, there exists $C_\varepsilon>0$
such that
\begin{equation}\label{eq:gaussian-absorption}
(1+t)^2\le C_\varepsilon e^{\varepsilon t^2},
\qquad t\ge0.
\end{equation}
  Set
\[
 \varphi_{0}=\frac1N\sum_{j=1}^N\varphi_j.
\]
Then the following identity holds:
\begin{equation}\label{eq:variance-identity}
 \sum_{j=1}^N|\varphi_j|^2
 =N|\varphi_{0}|^2
 +\frac1N\sum_{j<\ell}|\varphi_j-\varphi_\ell|^2.
\end{equation}
Applying \eqref{eq:gaussian-absorption} with
$\varepsilon=\alpha/(4N)$ to
\eqref{eq:vandermonde-growth}, and using
\eqref{eq:variance-identity}, we obtain 
\begin{equation} 
 |P(z)|
 \exp\!\left\{
 \frac{\alpha N}{2}|\varphi_{0}(z)|^2
 +\frac{\alpha}{4N}\sum_{j<\ell}|h_{j\ell}(z)|^2
 -\frac{\alpha N}{2}|z|^2
 \right\}
 \le C.  \nonumber 
\end{equation}
Fix $r<s$.  Discarding the other nonnegative terms in the exponent gives
\[
 |P(z)|
 e^{\frac{\alpha}{4N}|h_{rs}(z)|^2-\frac{\alpha N}{2}|z|^2}
 \le C.
\]
Since $r<s$ was arbitrary, the preceding
growth criterion yields
\[
h_{j\ell}(z)=a_{j\ell}z+b_{j\ell},
\qquad j\ne\ell,
\]
for some $a_{j\ell},b_{j\ell}\in\C$.
Since the symbols are pairwise distinct, each $h_{j\ell}$ is nonzero. If \(a_{j\ell}=0\),
then \(|h_{j\ell}(z)|=|b_{j\ell}|>0\).
Otherwise,
\[
|h_{j\ell}(z)|
\ge |a_{j\ell}||z|-|b_{j\ell}|
\longrightarrow\infty
\quad (|z|\to\infty).
\]

The family $\{h_{j\ell}:j\ne\ell\}$
consists of finitely many nonzero affine
functions. Hence, there exist
$R,\delta>0$ such that
\[
|h_{j\ell}(z)|\ge\delta
\qquad (|z|\ge R,\ j\ne\ell).
\]
Recalling that \(\rho(t)=t/(1+t)\) is
increasing, we obtain
\[
\prod_{\ell\ne j}\rho(|h_{j\ell}(z)|)
\ge
\left(\frac{\delta}{1+\delta}\right)^{N-1}
=:c>0
\]
for every \(j\) and \(|z|\ge R\).
Consequently, \eqref{eq:amplitude-with-separation}
gives
\[
m_{u_j,\varphi_j}(z)
\le \frac{C_{N,\alpha}}{c}\|T\|,
\qquad |z|\ge R.
\]
By continuity on the remaining compact disk,
we conclude that
\begin{equation} 
m_{u_j,\varphi_j}\in L^\infty(\mathbb C)
\qquad (1\le j\le N). \nonumber
\end{equation}

Since $u_j\not\equiv0$ and
\[
 |u_j(z)|e^{\frac{\alpha}{2}|\varphi_j(z)|^2-\frac{\alpha}{2}|z|^2}
 \le \norm{m_{u_j,\varphi_j}}_\infty,
\]
applying the same growth observation shows that every $\varphi_j$ is affine.

 Suppose first that
$\varphi_j(z)=a_jz+b_j$ with $a_j\ne0$.
Since $m_{u_j,\varphi_j}\in L^\infty(\C)$,
the change of variables $w=a_jz+b_j$ yields
\begin{align*}
 \norm{W_{u_j,\varphi_j}f}_\alpha^2
 &=\frac{\alpha}{\pi}\int_\C
 m_{u_j,\varphi_j}(z)^2|f(\varphi_j(z))|^2
 e^{-\alpha|\varphi_j(z)|^2}\,dA(z)\\
 &\le |a_j|^{-2}\norm{m_{u_j,\varphi_j}}_\infty^2\norm{f}_\alpha^2.
\end{align*}

It remains to consider the constant symbols.
Let
\[
I_0=\{j:\varphi_j(z)\equiv b_j
\text{ for some }b_j\in\mathbb C\}.
\]
If $I_0=\varnothing$, the proof is complete. 
Otherwise, since all terms with nonconstant composition symbols have
already been shown to be bounded, the operator
\[
S=T-\sum_{j\notin I_0}W_{u_j,\varphi_j}
\]
is bounded on $\F_\alpha^2$. Moreover,
\[
(Sf)(z)=\sum_{j\in I_0}u_j(z)f(b_j).
\]

Since the composition symbols are pairwise
distinct, the points $b_j$, $j\in I_0$, are
distinct. Thus, for each $j\in I_0$, there
exists a Lagrange interpolation polynomial
$q_j$ such that
\[
q_j(b_\ell)=
\begin{cases}
1,&\ell=j,\\
0,&\ell\ne j.
\end{cases}
\]

Applying $S$ to $q_j$, we obtain
\[
Sq_j
=\sum_{\ell\in I_0}u_\ell q_j(b_\ell)
=u_j.
\]
Since $q_j\in\F_\alpha^2$ and $S$ is bounded,
it follows that $u_j\in\F_\alpha^2$.

Finally, for every $f\in\F_\alpha^2$,
\[
W_{u_j,b_j}f=u_jf(b_j).
\]
Since point evaluation at $b_j$ is bounded
on $\F_\alpha^2$, we have
\[
\|W_{u_j,b_j}f\|_\alpha
\le e^{\frac{\alpha}{2}|b_j|^2}
\|u_j\|_\alpha\|f\|_\alpha.
\]
Hence, each term with constant composition symbol is a
bounded operator of rank at most one.
This completes the proof.
\end{proof}

\section{Compactness}
\label{sec:compactness}

For a bounded operator $A$ on $\F_\alpha^2$,
its essential norm is defined by
\begin{equation}
\|A\|_e
=\inf\{\|A-K\|:K\in\mathcal K(\F_\alpha^2)\}.
\nonumber
\end{equation}

We first recall the compactness criterion
for a single weighted composition operator.
If $\varphi$ is nonconstant, then
$W_{u,\varphi}$ is compact on $\F_\alpha^2$
if and only if
\[
m_{u,\varphi}(z)\longrightarrow0
\qquad (|z|\to\infty);
\]
see \cite[Theorem~2.3]{Le2014}.

Suppose further that $W_{u,\varphi}$
is bounded, $u\not\equiv0$, and
$\varphi$ is nonconstant. Then
$\varphi(z)=az+b$ with $0<|a|\le1$.
Taking $p=q=2$ in
\cite[Theorem~3.8]{TienKhoiWeighted2019},
we obtain
\begin{equation}\label{eq:single-essential} 
\limsup_{|z|\to\infty}m_{u,\varphi}(z)
 \le \|W_{u,\varphi}\|_e \le 2|a|^{-1}
\limsup_{|z|\to\infty}m_{u,\varphi}(z). 
\end{equation}

\begin{lemma}\label{lem:exterior-gram}
Let $M\ge1$ and let
$\phi_j(z)=a_jz+b_j$, where $a_j\ne0$,
$1\le j\le M$, be pairwise distinct
  functions.
Then there exist constants $R_0,c_0>0$
such that the Gram matrix
\[
G(z)=
\bigl(
\ip{k_{\phi_i(z)}}{k_{\phi_j(z)}}
\bigr)_{i,j=1}^M
\]
satisfies
\begin{equation}\label{eq:block-lower-general}
G(z)\ge c_0I_M
\qquad(|z|\ge R_0).
\end{equation}
\end{lemma}

\begin{proof}
Since the affine symbols are pairwise distinct,
for each $i\ne j$, the difference
\[
\phi_i(z)-\phi_j(z)
=(a_i-a_j)z+b_i-b_j
\]
is either a nonzero constant or tends to infinity
 as $|z|\to\infty$.
Since there are only finitely many pairs,
there exist $R_0,\delta>0$ such that
\[
|\phi_i(z)-\phi_j(z)|\ge\delta
\qquad(|z|\ge R_0,\ i\ne j).
\]

Fix $|z|\ge R_0$ and
$\xi=(\xi_1,\ldots,\xi_M)\in\C^M$.
Applying Lemma~\ref{lem:kernel-extraction}
to the normalized kernels, we obtain,
for every $1\le j\le M$,
\[
\left\|
\sum_{\ell=1}^M
\overline{\xi_\ell}k_{\phi_\ell(z)}
\right\|_\alpha
\ge C_{M,\alpha}^{-1}
\left(\frac{\delta}{1+\delta}\right)^{M-1}
|\xi_j|.
\]
Consequently,
\[
\begin{aligned}
\xi^*G(z)\xi
&=
\left\|
\sum_{j=1}^M
\overline{\xi_j}k_{\phi_j(z)}
\right\|_\alpha^2\\
&\ge
\frac{1}{M C_{M,\alpha}^2}
\left(\frac{\delta}{1+\delta}\right)^{2(M-1)}
\sum_{j=1}^M|\xi_j|^2.
\end{aligned}
\]
Thus $G(z)\ge c_0I_M$ for $|z|\ge R_0$,
where
\[
c_0=
\frac{1}{M C_{M,\alpha}^2}
\left(\frac{\delta}{1+\delta}\right)^{2(M-1)}
>0.
\]
\end{proof}

\begin{proof}[Proof of Theorem~\ref{T1.2}]
We first establish the essential norm
estimates.
Assume that $T$ is bounded.
By Theorem~\ref{T1.1}, every
$W_{u_j,\varphi_j}$ is bounded. 
If $u_j\equiv0$, then
$W_{u_j,\varphi_j}=0$ and
$m_{u_j,\varphi_j}\equiv0$. 
If $\varphi_j\equiv b_j$, then
$W_{u_j,b_j}$ has rank at most one.
Moreover,
\[
u_j=W_{u_j,b_j}1\in\F_\alpha^2,
\]
and hence
\[
\begin{aligned}
m_{u_j,b_j}(z)
&=e^{\alpha|b_j|^2/2}
|u_j(z)|e^{-\alpha|z|^2/2}\\
&=e^{\alpha|b_j|^2/2}
|\ip{u_j}{k_z}|
\longrightarrow0
\end{aligned}
\]
as $|z|\to\infty$, since
$k_z$ converges weakly to $0$ in $\F_\alpha^2$.

Thus, removing the zero and
constant symbol summands changes
neither the essential norm nor
the limit superior in
\eqref{eq:finite-essential}.    
If no summands remain, both the essential
norm and the limit superior vanish,
and the conclusion is immediate. 
Without loss of generality, we may assume that
all the remaining composition symbols are
nonconstant affine functions of the form
\[
\varphi_j(z)=a_jz+b_j,
\qquad a_j\ne0.
\]

By \eqref{eq:adjoint-kernel} and
Lemma~\ref{lem:exterior-gram},
there exist $R_0,c_0>0$ such that
\[
\|T^*k_z\|_\alpha^2
\ge c_0
\sum_{j=1}^N
m_{u_j,\varphi_j}(z)^2,
\qquad |z|\ge R_0.
\]

For every compact operator $K$,
we have
\[
\|K^*k_z\|_\alpha\longrightarrow0
\qquad(|z|\to\infty).
\]
Consequently,
\[
\begin{aligned}
\|T-K\|
&\ge
\limsup_{|z|\to\infty}
\|(T-K)^*k_z\|_\alpha\\
&\ge
\limsup_{|z|\to\infty}
\|T^*k_z\|_\alpha\\
&\ge
\sqrt{c_0}
\limsup_{|z|\to\infty}
\left(
\sum_{j=1}^N
m_{u_j,\varphi_j}(z)^2
\right)^{1/2}.
\end{aligned}
\]

Taking the infimum over compact
operators $K$ yields
\[
\|T\|_e
\ge
\sqrt{c_0}
\limsup_{|z|\to\infty}
\left(
\sum_{j=1}^N
m_{u_j,\varphi_j}(z)^2
\right)^{1/2}.
\]
This gives the lower bound estimate.

For the upper bound, set
\[
L=\limsup_{|z|\to\infty}
\left(
\sum_{j=1}^N m_{u_j,\varphi_j}(z)^2
\right)^{1/2}.
\]
Since
\[
\limsup_{|z|\to\infty}
m_{u_j,\varphi_j}(z)\le L
\qquad(1\le j\le N),
\]
the estimate \eqref{eq:single-essential} gives
\[
\begin{aligned}
\|T\|_e
&\le \sum_{j=1}^N
\|W_{u_j,\varphi_j}\|_e\\
&\le 2\sum_{j=1}^N |a_j|^{-1}
\limsup_{|z|\to\infty}
m_{u_j,\varphi_j}(z)\\
&\le 2\left(
\sum_{j=1}^N |a_j|^{-1}
\right)L.
\end{aligned}
\]

The compactness assertion follows immediately
from \eqref{eq:finite-essential} and the compactness criterion for individual weighted composition operators. 
This completes the proof.
\end{proof}

\section{Schatten-class membership}
\label{sec:finite-schatten}

We recall the Schatten-class criteria
for individual weighted composition
operators on Fock spaces
\cite{Mengestie2013,Mengestie2015}.
Let $u,\varphi\in\mathcal H(\C)$, with
$\varphi$ nonconstant. For every
$0<p<\infty$,
\[
W_{u,\varphi}\in\Sp_p
\iff m_{u,\varphi}\in L^p(\C,dA).
\]
Moreover, when these conditions hold,
\begin{equation}\label{eq:single-schatten-upper}
\|W_{u,\varphi}\|_{\Sp_p}^p
\le C\int_\C m_{u,\varphi}(z)^p\,dA(z),
\end{equation}
where $C$ depends only on $\alpha$,
$p$, and the fixed symbol $\varphi$.

For $0<p\le2$, this estimate follows from
the reproducing kernel trace formula and
the concavity of $t\mapsto t^{p/2}$.
For $p\ge2$, it follows by representing
$W_{u,\varphi}^*W_{u,\varphi}$ as a positive
Toeplitz operator and applying the standard
Schatten-class estimate for such operators.
In both cases, the constant $C$ is
independent of the weight $u$.

\begin{lemma}\label{lem:kernel-family-estimates}
Let $M\ge1$, and 
\[
\phi_j(z)=a_jz+b_j,\qquad a_j\ne0,
\quad 1\le j\le M,
\]
be pairwise distinct affine maps. Choose
$R_0,c_0>0$ as in Lemma~\ref{lem:exterior-gram}.

\begin{enumerate}[label=\textup{(\roman*)}]
\item For every relatively separated set
$Z\subset\C$ and every $R>0$, the set
$Z\setminus D(0,R_0)$ can be partitioned
into finitely many classes such that,
for distinct $z,w$ in the same class,
\[
|z-w|\ge R,
\qquad
|\phi_i(z)-\phi_j(w)|\ge R
\quad(1\le i,j\le M).
\]

\item For sufficiently large $R$,
suppose that $\{z_k\}_{k\in J}$
satisfies $|z_k|\ge R_0$ and
\[
|\phi_i(z_k)-\phi_j(z_\ell)|\ge R
\qquad
(k\ne\ell,\ 1\le i,j\le M).
\]
Then, for every finitely supported family
of complex coefficients $(c_{j,k})$,
\begin{equation}\label{eq:kernel-family-estimate}
\begin{aligned}
\frac{c_0}{2}
\sum_{k\in J}\sum_{j=1}^M|c_{j,k}|^2
&\le
\left\|
\sum_{k\in J}\sum_{j=1}^M
c_{j,k}k_{\phi_j(z_k)}
\right\|_\alpha^2 \le
\left(M+\frac{c_0}{2}\right)
\sum_{k\in J}\sum_{j=1}^M|c_{j,k}|^2.
\end{aligned}
\end{equation}
 
\end{enumerate}
\end{lemma}

\begin{proof}

\emph{(i)}
Fix $R>0$ and put
\[
Z_0=Z\setminus D(0,R_0).
\]

Since $Z$ is relatively separated,
for every $r>0$ there exists an integer
$N_r<\infty$ such that
\begin{equation}\label{eq:relative-separation-count}
\sup_{\zeta\in\C}
\#(Z\cap D(\zeta,r))\le N_r.
\end{equation}
Indeed, every disk of fixed radius
can be covered by finitely many unit disks.

Fix $z\in Z_0$. We estimate the number
of points $w\in Z_0\setminus\{z\}$
for which at least one of the inequalities
\[
|z-w|\ge R,
\qquad
|\phi_i(z)-\phi_j(w)|\ge R
\quad(1\le i,j\le M)
\]
fails.

If $|z-w|<R$, then $w\in D(z,R)$. 
If $|\phi_i(z)-\phi_j(w)|<R$ for some $i,j$, since
$\phi_j(w)=a_jw+b_j$ and $a_j\ne0$,
we obtain
\[
w\in D\left(
\frac{a_iz+b_i-b_j}{a_j},
\frac{R}{|a_j|}
\right).
\]
Consequently, by
\eqref{eq:relative-separation-count},
the number of such points $w$
is bounded above by
\[
K:=N_R+\sum_{i,j=1}^M
N_{R/|a_j|}<\infty,
\]
which is independent of $z$.

For distinct $z,w\in Z_0$, write $z\sim w$
if $|z-w|<R$
 or $
|\phi_i(z)-\phi_j(w)|<R$,
 for some $1\le i,j\le M$.
This relation is symmetric, and the preceding
counting estimate shows that each point of
$Z_0$ is related to at most $K$ other points.

We construct the desired partition by
successive maximal selections.
Set $Y_1=Z_0$.
Since $Z_0$ is at most countable, we may
choose a maximal subset $Z_1\subset Y_1$
such that no two distinct points of $Z_1$
are related.
Indeed, such a subset can be obtained by
enumerating $Y_1$ and successively retaining
each point that is not related to any
previously retained point.

Set $Y_2=Y_1\setminus Z_1$.
Inductively, having constructed $Y_s$,
choose a maximal subset $Z_s\subset Y_s$
whose distinct points are not related,
and put
\[
Y_{s+1}=Y_s\setminus Z_s.
\]

We claim that $Y_{K+2}=\varnothing$.
Otherwise, choose $z\in Y_{K+2}$.
For each $1\le s\le K+1$, we have
$z\in Y_s\setminus Z_s$.
By the maximality of $Z_s$, there exists
$w_s\in Z_s$ such that $z\sim w_s$. Since the sets $Z_s$ are pairwise disjoint,
the points $w_1,\ldots,w_{K+1}$ are distinct.
Thus $z$ is related to at least $K+1$
distinct points, contradicting the
uniform bound $K$.

Consequently,
\[
Z_0=Z_1\cup\cdots\cup Z_L,
\qquad L\le K+1,
\]
where empty sets are omitted.
By construction, any two distinct points
$z,w$ in the same class satisfy
\[
|z-w|\ge R,
\qquad
|\phi_i(z)-\phi_j(w)|\ge R
\quad(1\le i,j\le M).
\]
This proves \textup{(i)}.

\emph{(ii)}
We consider the Gram matrix $\mathcal G$
of the family $\{k_{\phi_j(z_k)}:
1\le j\le M,\ k\in J\}$, 
acting on
$\ell^2(\{1,\ldots,M\}\times J)$.
Its entries are given by
\[
\mathcal G_{(i,k),(j,\ell)}
=
\langle k_{\phi_i(z_k)},
k_{\phi_j(z_\ell)}\rangle_\alpha.
\]

Write $\mathcal G=B+E$, where
$E=\mathcal G-B$ and $B$ is defined entrywise by
\[
B_{(i,k),(j,\ell)}
=
\begin{cases}
\langle k_{\phi_i(z_k)},
k_{\phi_j(z_k)}\rangle_\alpha,
& k=\ell,\\
0, & k\ne\ell.
\end{cases}
\]

Thus $B$ is block diagonal, with
\[
B=\bigoplus_{k\in J}B_k,
\qquad
B_k=
\bigl(
\langle k_{\phi_i(z_k)},
k_{\phi_j(z_k)}\rangle_\alpha
\bigr)_{i,j=1}^M.
\]

Since $|z_k|\ge R_0$,
Lemma~\ref{lem:exterior-gram} gives
\[
B_k\ge c_0I_M.
\]
Moreover, $B_k$ is positive semidefinite
and its diagonal entries are all equal
to one. Therefore,
\[
\|B_k\|
\le\operatorname{tr}(B_k)=M.
\]
It follows that
\begin{equation}\label{eq:gram-block-bounds}
c_0I\le B\le MI.
\end{equation}

We next estimate the off-diagonal part $E$.
Fix $i,k$ and $j$.
By the separation hypothesis,
the points $
\{\phi_j(z_\ell):\ell\ne k\}$ 
are $R$-separated. Indeed, for
$\ell\ne\ell'$,
\[
|\phi_j(z_\ell)-\phi_j(z_{\ell'})|
\ge R.
\] 

For $n\ge1$, consider the annulus
\[
A_n=
\{w\in\C:
nR\le|w-\phi_i(z_k)|<(n+1)R\}.
\]
Since the points are $R$-separated,
the disks of radius $R/2$ centered
at these points are pairwise disjoint.

Those corresponding to points in $A_n$
are contained in $
D\bigl(\phi_i(z_k),(n+3/2)R\bigr)$. 
Comparing areas, we find that
\[
\#\{\ell\ne k:\phi_j(z_\ell)\in A_n\}
\le C(n+1)^2,
\]
where $C>0$ is an absolute constant.

Consequently,
\[ 
 \sum_{\ell\ne k}
e^{-\alpha|\phi_i(z_k)-\phi_j(z_\ell)|^2/2} \le
C\sum_{n=1}^\infty
(n+1)^2e^{-\alpha n^2R^2/2} =:\varepsilon_\alpha(R). 
\]
The Gaussian decay implies that
\[
\varepsilon_\alpha(R)\longrightarrow0
\qquad(R\to\infty).
\]

Summing over $j=1,\ldots,M$ and using
\eqref{eq:kernel-overlap}, we obtain
\[
\sup_{i,k}
\sum_{\ell\ne k}\sum_{j=1}^M
\big|
\langle k_{\phi_i(z_k)},
k_{\phi_j(z_\ell)}\rangle_\alpha
\big|
\le M\varepsilon_\alpha(R).
\]
Since the Gram matrix is Hermitian,
the same estimate holds for its column
sums. The Schur test therefore yields
\[
\|E\|\le M\varepsilon_\alpha(R).
\]

Choose $R$ sufficiently large such that
\[
M\varepsilon_\alpha(R)\le\frac{c_0}{2}.
\]
Combining this with
\eqref{eq:gram-block-bounds}
gives
\[
\frac{c_0}{2}I
\le\mathcal G
\le\left(M+\frac{c_0}{2}\right)I.
\]

Thus, for every finitely supported family
$(c_{j,k})$, we have
\[ 
\frac{c_0}{2}
\sum_{k\in J}\sum_{j=1}^M|c_{j,k}|^2
 \le
\left\|
\sum_{k\in J}\sum_{j=1}^M
c_{j,k}k_{\phi_j(z_k)}
\right\|_\alpha^2 \le
\left(M+\frac{c_0}{2}\right)
\sum_{k\in J}\sum_{j=1}^M|c_{j,k}|^2. 
\]
This proves
\eqref{eq:kernel-family-estimate}.
 
\end{proof}

\begin{proof}[Proof of Theorem \ref{T1.3}]
If all the summands belong to $\Sp_p$, then so does their
finite sum, including when $0<p<1$. We prove the converse
and the quantitative lower estimate.

Suppose that $T\in\Sp_p$. By Theorem~\ref{T1.1},
every summand is bounded.
Omit all zero summands and subtract the summands
with constant composition symbols.
The removed summands have finite rank, so the
remaining sum still belongs to $\Sp_p$.
If no summands remain, the conclusion is immediate.

Relabel the remaining summands and continue to
denote their number by $N$ and their sum by $T$.
By the boundedness criterion for individual
weighted composition operators, we have
$u_j\not\equiv0$ and
\[
\varphi_j(z)=a_jz+b_j,\qquad a_j\ne0,
\quad 1\le j\le N.
\]

Take
\[
F(z)=\sum_{j=1}^Nm_{u_j,\varphi_j}(z)^p.
\]
Let $\{Q_\nu\}_{\nu\in\mathbb Z^2}$
be a partition of $\C$ into half-open unit squares  and choose
$z_\nu\in Q_\nu$ such that
\begin{equation}\label{eq:finite-max-sampling}
F(z_\nu)\ge\frac12\sup_{z\in Q_\nu}F(z).
\end{equation}
Since each $Q_\nu$ has area one, this choice implies
\[
\int_{Q_\nu}F(z)\,dA(z)
\le\sup_{z\in Q_\nu}F(z)
\le2F(z_\nu).
\]
Moreover, every unit disk meets only a uniformly
bounded number of unit squares. Thus
$Z=\{z_\nu\}$ is relatively separated, with a
density bound independent of the weights.

Choose $R_0,c_0>0$ from Lemma~\ref{lem:exterior-gram}
for the fixed family $\{\varphi_j:1\le j\le N\}$.
Fix $R$ sufficiently large for
Lemma~\ref{lem:kernel-family-estimates}(ii) to hold.
By part~(i) of that lemma, the points of
$Z\setminus D(0,R_0)$ can be partitioned into
finitely many classes satisfying its separation
hypotheses. Since $Z$ is relatively separated,
the number of classes and the constants in 
\eqref{eq:kernel-family-estimate} depend only on
$\alpha$ and the fixed symbols, not on the weights.

Fix one nonempty class $Z_r=\{z_k:k\in J\}$.
By \eqref{eq:adjoint-kernel},
\begin{equation}\label{eq:finite-Tstar}
T^*k_{z_k}
=\sum_{j=1}^N\xi_{j,k}k_{\varphi_j(z_k)},
\qquad
|\xi_{j,k}|=m_{u_j,\varphi_j}(z_k).
\end{equation}
Put
\[
d_k=\left(\sum_{j=1}^N|\xi_{j,k}|^2\right)^{1/2}.
\]
Indices with $d_k=0$ contribute nothing to $F(z_k)$
and may be omitted. For all remaining indices, define
\[
h_k=\frac1{d_k}\sum_{j=1}^N
\xi_{j,k}k_{\varphi_j(z_k)}.
\]
Then $T^*k_{z_k}=d_kh_k$. For every finitely
supported sequence $(t_k)$, we have
\[
\begin{aligned}
\sum_{k\in J}\sum_{j=1}^N
\left|t_k\frac{\xi_{j,k}}{d_k}\right|^2
&=\sum_{k\in J}\frac{|t_k|^2}{d_k^2}
\sum_{j=1}^N|\xi_{j,k}|^2\\
&=\sum_{k\in J}|t_k|^2.
\end{aligned}
\]
Applying \eqref{eq:kernel-family-estimate} with
$c_{j,k}=t_k\xi_{j,k}/d_k$ therefore gives
\[
\frac{c_0}{2}\sum_{k\in J}|t_k|^2
\le
\left\|\sum_{k\in J}t_kh_k\right\|_\alpha^2
\le\left(N+\frac{c_0}{2}\right)
\sum_{k\in J}|t_k|^2.
\]
The sample points $(z_k)$ are $R$-separated.
Hence, the estimate for separated normalized
reproducing kernels in
Section~\ref{sec:preliminaries} also yields 
\[
\left\|\sum_{k\in J}t_kk_{z_k}\right\|_\alpha^2
\le C_{\alpha,R}\sum_{k\in J}|t_k|^2.
\]

Define operators $U,V:\ell^2(J)\to\F_\alpha^2$ by
\[
Ue_k=k_{z_k},\qquad Ve_k=h_k,
\]
and let $D$ be the diagonal operator $De_k=d_ke_k$.
The two estimates above show that $U,V$ are bounded
and
\[
V^*V\ge\frac{c_0}{2}I.
\]
Consequently,
\[
L_V=(V^*V)^{-1}V^*
\]
is a bounded left inverse of $V$, satisfying
\[
L_VV=I,\qquad \|L_V\|\le\sqrt{\frac2{c_0}}.
\]
Furthermore, by Lemma~\ref{lem:exterior-gram} and
\eqref{eq:finite-Tstar},
\[
c_0d_k^2\le\|T^*k_{z_k}\|_\alpha^2
\le\|T\|^2.
\]
It follows that $D$ is bounded. The identity
$T^*Ue_k=d_kh_k=VDe_k$ holds on every coordinate
vector. Since all the operators involved are bounded,
it extends from finitely supported sequences to
$\ell^2(J)$. Applying the left inverse of $V$ yields
\begin{equation}\label{eq:discrete-factorization}
T^*U=VD,\qquad D=L_VT^*U.
\end{equation}
Since $T^*\in\Sp_p$, the two-sided ideal property
\eqref{eq:ideal-property} shows that $D\in\Sp_p$.
As $D$ is nonnegative and diagonal, its singular
values are precisely the numbers $d_k$, counted with
multiplicity. Hence,
\begin{equation}\label{eq:finite-d-lp}
\sum_{k\in J}d_k^p
=\|D\|_{\Sp_p}^p
\le \|L_V\|^p\|U\|^p\|T\|_{\Sp_p}^p.
\end{equation}
This remains valid for $0<p<1$; no triangle
inequality for the Schatten quasi-norm is used here.

An elementary comparison of finite sums gives
\[
\begin{aligned}
F(z_k)
&=\sum_{j=1}^N|\xi_{j,k}|^p\\
&\le N^{\max\{1-p/2,0\}}
\left(\sum_{j=1}^N|\xi_{j,k}|^2\right)^{p/2}\\
&=N^{\max\{1-p/2,0\}}d_k^p.
\end{aligned}
\]
By \eqref{eq:finite-d-lp}, the sum of $F(z_k)$
over each class is at most a fixed multiple of
$\|T\|_{\Sp_p}^p$. There are only finitely many
classes, with a uniformly bounded number of classes.
Combining these estimates with
\eqref{eq:finite-max-sampling}, we obtain
\begin{equation}\label{eq:finite-exterior-integral}
\sum_{\nu:\,|z_\nu|\ge R_0}
\int_{Q_\nu}F(z)\,dA(z)
\le C\|T\|_{\Sp_p}^p.
\end{equation}
Here and below, the constants depend only on
$\alpha,p$ and the fixed composition symbols.

It remains to bound the contribution of the finitely
many squares whose sample points satisfy
$|z_\nu|<R_0$. Each such square is contained in
$D(0,R_0+\sqrt2)$. Fix $R_1>R_0+\sqrt2$.
For $|w|=R_1$, Lemma~\ref{lem:exterior-gram}
and \eqref{eq:adjoint-kernel} give
\[
c_0\sum_{j=1}^Nm_{u_j,\varphi_j}(w)^2
\le \|T^*k_w\|_\alpha^2
\le \|T\|^2.
\]
In particular,
\[
m_{u_j,\varphi_j}(w)
\le c_0^{-1/2}\|T\|,\qquad 1\le j\le N.
\]
Consequently,
\[
\sup_{|w|=R_1}|u_j(w)|
\le c_0^{-1/2}\|T\|
\sup_{|w|=R_1}
 e^{\frac\alpha2(|w|^2-|\varphi_j(w)|^2)}.
\]
The maximum-modulus principle gives
$\sup_{|z|\le R_1}|u_j(z)|
\le\sup_{|w|=R_1}|u_j(w)|$.
Multiplying by the fixed Gaussian factor in the
definition of $m_{u_j,\varphi_j}$, we obtain
\begin{equation} 
\sup_{|z|\le R_1}m_{u_j,\varphi_j}(z)
\le C_j\|T\|
\le C_j\|T\|_{\Sp_p},
\qquad 1\le j\le N,  \nonumber
\end{equation}
where $C_j$ depends only on $\alpha,R_1,c_0$ and
$\varphi_j$, and not on $u_j$.
Since the omitted squares are disjoint and contained
in $D(0,R_1)$, it follows that
\[
\sum_{\nu:\,|z_\nu|<R_0}
\int_{Q_\nu}F(z)\,dA(z)
\le \pi R_1^2\sum_{j=1}^NC_j^p
\|T\|_{\Sp_p}^p.
\]
Together with \eqref{eq:finite-exterior-integral},
this proves
\begin{equation} 
\sum_{j=1}^N\int_\C
m_{u_j,\varphi_j}(z)^p\,dA(z)
\le C\|T\|_{\Sp_p}^p.
    \nonumber
\end{equation}
In particular, each $m_{u_j,\varphi_j}$ belongs to
$L^p(\C,dA)$. The estimate for individual operators in \eqref{eq:single-schatten-upper} now implies 
\[
\sum_{j=1}^N
\|W_{u_j,\varphi_j}\|_{\Sp_p}^p
\le C'\|T\|_{\Sp_p}^p.
\]
The discarded summands also belong to $\Sp_p$.
Together with the preceding argument for the
remaining summands, this proves the equivalence
for the original sum.

If all the original composition symbols are
nonconstant, only zero summands have been omitted,
so the operator $T$ has not changed.
The preceding estimate therefore gives the lower
bound in \eqref{eq:finite-schatten-quantitative}.
The upper bound follows from the finite-sum
inequality in Section~\ref{sec:preliminaries}.
This completes the proof.

\end{proof}

By the proof of Theorem \ref{T1.3} and the Schatten class characterization for the weighted composition operator, we have the following corollary.
\begin{corollary}\label{cor:intro-amplitude}
Let $0<p<\infty$ and $T$ be given by
\eqref{eq:formal-sum}, where
$\varphi_1,\ldots,\varphi_N$ are pairwise
distinct nonconstant entire functions.
Then $T\in\Sp_p$ if and only if $m_{u_j,\varphi_j}\in L^p(\C,dA),$ for $1\le j\le N.$
\end{corollary}

\section{Finite sums of generalized Volterra-type integral operators}
\label{sec:volterra}

We consider finite sums of generalized Volterra-type integral operators
\[
V_{(g,\phi)}f(z)
=\int_0^z f(\phi(\zeta))g'(\zeta)\,d\zeta,
\qquad g,\phi\in \mathcal H(\C),
\]
where the integral is path independent.  Put
\begin{equation} 
\mathcal V=\sum_{j=1}^N V_{(g_j,\phi_j)}.
    \nonumber 
\end{equation} 

Recall that
\[
f'(w)=\langle f,K_w^{[1]}\rangle_\alpha,
\qquad K_w^{[1]}(z)=\alpha z e^{\alpha z\overline w},
\]
and
\begin{equation}\label{e1}
\|K_w^{[1]}\|_\alpha^2
=\alpha(1+\alpha|w|^2)e^{\alpha|w|^2}.
\end{equation}
Set \(\ell_w=K_w^{[1]}/\|K_w^{[1]}\|_\alpha\). Then
\begin{equation}\label{eq:derivative-kernel-overlap}
\big|\langle\ell_z,\ell_w\rangle_\alpha\big|
=\frac{|1+\alpha w\overline z|}
{\sqrt{(1+\alpha|z|^2)(1+\alpha|w|^2)}}
 e^{-\alpha|z-w|^2/2}
\le e^{-\alpha|z-w|^2/2}.
\end{equation}
Consequently, for every $\delta>0$ there exists
a constant $C_{\alpha,\delta}>0$ such that
\[
\left\|\sum_k c_k\ell_{w_k}\right\|_\alpha^2
\le C_{\alpha,\delta}\sum_k|c_k|^2
\]
for every finitely supported $(c_k)$,
whenever the points $(w_k)$ are
$\delta$-separated. 
Moreover, $\ell_w$ converges weakly to $0$
in $\F_\alpha^2$ as $|w|\to\infty$.
This follows first for polynomials by
derivative evaluation and then for all
functions in $\F_\alpha^2$ by density.

For entire functions $g$ and $\phi$, set
\begin{equation}\label{eq:volterra-amplitude}
\mu_{g,\phi}(w)
=
\frac{|g'(w)|}
{\sqrt{\alpha(1+\alpha|w|^2)}}
\exp\left(
\frac{\alpha}{2}
(|\phi(w)|^2-|w|^2)
\right).
\end{equation}

Whenever $V_{(g,\phi)}$ is bounded,
differentiation and the reproducing
property give
\begin{equation}\label{eq:volterra-adjoint-derivative}
V_{(g,\phi)}^*\ell_w
=
\frac{\overline{g'(w)}}
{\sqrt{\alpha(1+\alpha|w|^2)}}
\exp\left(
\frac{\alpha}{2}
(|\phi(w)|^2-|w|^2)
\right)k_{\phi(w)}.
\end{equation} 
In particular, whenever $V_{(g,\phi)}$
is bounded,
\[
\|V_{(g,\phi)}^*\ell_w\|_\alpha
=\mu_{g,\phi}(w).
\]

If  \(\mathcal V\) is bounded, its adjoint
  satisfies the identity
\begin{equation}\label{eq:volterra-finite-adjoint}
\mathcal V^*K_w^{[1]}
=\sum_{j=1}^N\overline{g_j'(w)}K_{\phi_j(w)}.
\end{equation}

We recall the following results for individual
generalized Volterra-type integral operators
\cite{MengestieWorku2018}.
Let $g$ and $\phi$ be nonconstant entire functions.
Then
\[
\begin{aligned}
V_{(g,\phi)}\text{ is bounded on }\F_\alpha^2
&\iff \mu_{g,\phi}\in L^\infty(\C),\\
V_{(g,\phi)}\text{ is compact on }\F_\alpha^2
&\iff \mu_{g,\phi}(z)\to0
\quad (|z|\to\infty).
\end{aligned}
\]
In particular, either condition forces $\phi$
to be affine.
Moreover, if $V_{(g,\phi)}$ is bounded,
then, for every $0<p<\infty$,
\begin{equation}\label{eq:volterra-single-criterion}
V_{(g,\phi)}\in\Sp_p
\iff \mu_{g,\phi}\in L^p(\C,dA).
\end{equation}
For a constant symbol \(\phi\equiv b\), we use instead
\begin{equation} 
V_{(g,b)}f=f(b)(g-g(0)),\nonumber
\end{equation}
which has rank at most one whenever it is bounded.

\begin{proof}[Proof of Theorem~\ref{T1.4}]
Sufficiency is immediate for each of the
three properties. We prove necessity separately.

\emph{Boundedness.}
Suppose that $\mathcal V$ is bounded.
If all summands are zero, there is nothing to prove.
After omitting zero summands and relabeling,
we may assume that $g_j'\not\equiv0$
for $1\le j\le N$.
The case $N=1$ is immediate. By
\eqref{eq:volterra-finite-adjoint}, \eqref{e1}, and
Lemma~\ref{lem:kernel-extraction},
\begin{equation}\label{eq:volterra-coefficient-extraction}
\mu_{g_j,\phi_j}(z)
\prod_{\ell\ne j}
\frac{|\phi_j(z)-\phi_\ell(z)|}
{1+|\phi_j(z)-\phi_\ell(z)|}
\le C\|\mathcal V\|,
\qquad 1\le j\le N,
\end{equation}
where $C$ depends only on $\alpha$ and $N$.

Multiplying these inequalities and arguing as in the proof of
Theorem~\ref{T1.1}, we use \eqref{eq:variance-identity},
\eqref{eq:gaussian-absorption}, and
\[
[\alpha(1+\alpha|z|^2)]^{N/2}
\le C e^{\alpha N|z|^2/4}
\]
to obtain
\begin{equation} 
\left|\prod_{j=1}^N g_j'(z)
\prod_{j<\ell}(\phi_j(z)-\phi_\ell(z))^2\right|\exp\!\left(
\frac{\alpha}{4N}\sum_{j<\ell}|\phi_j(z)-\phi_\ell(z)|^2
-\frac{3\alpha N}{4}|z|^2\right)\le C.\nonumber
\end{equation}
For any $r<s$, the same inequality holds with the sum in the
exponent replaced by $|\phi_r(z)-\phi_s(z)|^2$.
The product outside the exponential is a nonzero entire function.
The growth estimate used in the proof of Theorem~\ref{T1.1}
therefore shows that every $\phi_r-\phi_s$ is affine.
Since the composition symbols are pairwise distinct, there are
$R,\delta>0$ such that
\[
|\phi_j(z)-\phi_\ell(z)|\ge\delta
\qquad(|z|\ge R,\ j\ne\ell).
\]
It follows from \eqref{eq:volterra-coefficient-extraction}
and continuity on the remaining disk that
\[
\mu_{g_j,\phi_j}\in L^\infty(\C),
\qquad 1\le j\le N.
\]
For each nonconstant $\phi_j$, an argument analogous to that in
\cite[Theorem~2.1(i)]{MengestieWorku2018} 
implies that $V_{(g_j,\phi_j)}$ is bounded,
since
$\sqrt{\alpha(1+\alpha|z|^2)}
\asymp_\alpha 1+|z|$. 
Moreover, the argument of
\cite[Proposition~5.1(i)]{MengestieWorku2018}
shows that each nonconstant $\phi_j$ is affine.

Finally, the remaining sum is a bounded
finite sum of weighted composition operators
with pairwise distinct constant symbols.
By Theorem~\ref{T1.1}, each summand is bounded
and hence has rank at most one.
This proves the boundedness assertion.

\emph{Compactness.}
Suppose that $\mathcal V$ is compact.
By the boundedness assertion just proved,
every summand is bounded. 
Omit all zero summands and subtract the summands
with constant composition symbols.
The removed summands have finite rank and are
therefore compact, and the remaining sum is compact.
If no summands remain, the conclusion is immediate.

Relabel the remaining summands and continue to
denote their number by $N$ and their sum by $\mathcal V$.
We then have $g_j'\not\equiv0$ for every $j$.
By the individual boundedness criterion, the
remaining composition symbols are pairwise distinct
nonconstant affine functions.
By \eqref{eq:volterra-adjoint-derivative} and
Lemma~\ref{lem:exterior-gram}, for all sufficiently large $|z|$,
\[
\|\mathcal V^*\ell_z\|_\alpha^2
\ge c_0\sum_{j=1}^N\mu_{g_j,\phi_j}(z)^2.
\] 
Since $\ell_z$ converges weakly to $0$
in $\F_\alpha^2$ as $|z|\to\infty$
and $\mathcal V^*$ is compact,
the left-hand side tends to zero as $|z|\to\infty$.
Hence, every $\mu_{g_j,\phi_j}(z)\to0$, and the
compactness criterion for individual operators  proves the compactness assertion.

\emph{Schatten-class membership.}
Fix $0<p<\infty$, and suppose that $\mathcal V\in\Sp_p$. 
By the boundedness assertion just proved,
every summand is bounded. 
Omit all zero summands and subtract the summands
with constant composition symbols.
The removed summands have finite rank and hence
belong to $\Sp_p$, so the remaining sum still
belongs to $\Sp_p$.
If no summands remain, the conclusion is immediate.

Relabel the remaining summands and continue to
denote their number by $N$ and their sum by $\mathcal V$.
We then have $g_j'\not\equiv0$ for every $j$.
By the individual boundedness criterion, the
remaining composition symbols are pairwise distinct
nonconstant affine functions. Let
\[
F(z)=\sum_{j=1}^N\mu_{g_j,\phi_j}(z)^p.
\]
As in the proof of Theorem~\ref{T1.3}, let
$\{Q_\nu\}$ be a partition of $\C$ into half-open unit
squares and choose $z_\nu\in Q_\nu$ such that
\[
F(z_\nu)\ge\frac12\sup_{z\in Q_\nu}F(z).
\]
Consequently, $\int_{Q_\nu}F\,dA\le2F(z_\nu)$.
By Lemma~\ref{lem:kernel-family-estimates}, the sample points
outside $D(0,R_0)$, with $R_0$ from
Lemma~\ref{lem:exterior-gram}, can be partitioned into
finitely many classes satisfying
\eqref{eq:kernel-family-estimate}.

Fix one such class $\{z_k:k\in J\}$. Then
\[
\mathcal V^*\ell_{z_k}
=\sum_{j=1}^N\xi_{j,k}k_{\phi_j(z_k)}
=d_kh_k,
\qquad |\xi_{j,k}|=\mu_{g_j,\phi_j}(z_k),
\]
where, for $d_k>0$,
\[
d_k=\left(\sum_{j=1}^N|\xi_{j,k}|^2\right)^{1/2},
\qquad
h_k=d_k^{-1}\sum_{j=1}^N\xi_{j,k}k_{\phi_j(z_k)}.
\]
Indices with $d_k=0$ may be omitted.
For finitely supported $(t_k)$,
\[
\sum_{k\in J}\sum_{j=1}^N
\left|t_k\frac{\xi_{j,k}}{d_k}\right|^2
=\sum_{k\in J}|t_k|^2.
\]
It follows from Lemma~\ref{lem:kernel-family-estimates} that
\[
\frac{c_0}{2}\sum_{k\in J}|t_k|^2
\le\left\|\sum_{k\in J}t_kh_k\right\|_\alpha^2
\le\left(N+\frac{c_0}{2}\right)\sum_{k\in J}|t_k|^2.
\]
Moreover, \eqref{eq:derivative-kernel-overlap} and the
separation of $(z_k)$ give
\[
\left\|\sum_{k\in J}t_k\ell_{z_k}\right\|_\alpha^2
\le C\sum_{k\in J}|t_k|^2.
\]
Define
\[
Ue_k=\ell_{z_k},\qquad Ve_k=h_k,\qquad De_k=d_ke_k.
\]
These estimates show that $U$ and $V$ are bounded,
while $V^*V\ge(c_0/2)I$. Thus
$L_V=(V^*V)^{-1}V^*$ is a bounded left inverse of $V$. 
Lemma~\ref{lem:exterior-gram} also implies that
$d_k\le c_0^{-1/2}\|\mathcal V\|$, so $D$ is bounded. 
As in \eqref{eq:discrete-factorization}, we have
\[
\mathcal V^*U=VD,\qquad D=L_V\mathcal V^*U.
\]
The two-sided ideal property \eqref{eq:ideal-property} yields
\[
\sum_{k\in J}d_k^p
=\|D\|_{\Sp_p}^p
\le\|L_V\|^p\|U\|^p\|\mathcal V\|_{\Sp_p}^p.
\]
Since
\[
F(z_k)=\sum_{j=1}^N|\xi_{j,k}|^p
\le N^{\max\{1-p/2,0\}}d_k^p,
\]
summing over the finitely many classes and using the
sampling estimate gives
\[
\sum_{\nu:\,|z_\nu|\ge R_0}
\int_{Q_\nu}F(z)\,dA(z)<\infty.
\]
Only finitely many squares remain, and their integrals
are finite because $F$ is continuous. Consequently,
\[
\mu_{g_j,\phi_j}\in L^p(\C,dA),
\qquad 1\le j\le N.
\]
The criterion for individual operators  in
\eqref{eq:volterra-single-criterion} then implies
$V_{(g_j,\phi_j)}\in\Sp_p$, for every $j$.
This completes the proof.
\end{proof}

\subsection*{Author contributions}

All authors contributed equally to the research,
development, and writing of this manuscript.
All authors read and approved the final manuscript.

\subsection*{Conflict of interest}

The authors declare that they have no conflicts
of interest.

\subsection*{Data availability statement}

No datasets were generated or analyzed during
the current study.

\subsection*{Funding}

Y. Lu was supported by the National Natural
Science Foundation of China (Grant Nos. 12031002
and 12671149). C. Zu was supported by the National
Natural Science Foundation of China (Grant No. 12401151).

\subsection*{AI disclosure}
During the preparation of this manuscript, the authors used ChatGPT (OpenAI) for English language editing and assistance in reviewing mathematical arguments, including identifying potential gaps and discussing possible proof strategies. All mathematical statements and proofs were independently verified by the authors, who take full responsibility for the content of the manuscript.

\bibliographystyle{amsplain}
\bibliography{references}

\medskip

\noindent
School of Mathematical Sciences, Dalian University of Technology,
Dalian, Liaoning 116024, P. R. China

\noindent
Email address: \texttt{linzhaopeng2606@163.com}

\medskip

\noindent
School of Mathematical Sciences, Dalian University of Technology,
Dalian, Liaoning 116024, P. R. China

\noindent
Email address: \texttt{YRong@mail.dlut.edu.cn}

\medskip

\noindent
School of Mathematical Sciences, Dalian University of Technology,
Dalian, Liaoning 116024, P. R. China

\noindent
Email address: \texttt{zuchao@dlut.edu.cn}

\medskip

\noindent
School of Mathematical Sciences, Dalian University of Technology,
Dalian, Liaoning 116024, P. R. China

\noindent
Email address: \texttt{lyfdlut@dlut.edu.cn}

\end{document}